\documentclass[10pt]{amsart}

\usepackage{amsfonts,amssymb,amsmath}
\usepackage{mathrsfs}
\usepackage[top=1in, bottom=1in, left=1.25in, right=1.25in]{geometry}
\usepackage{hyperref}
\usepackage{xypic}

\usepackage{enumerate}

\newtheorem{Thm}{Theorem}

\newcommand{\Z}{\mathbb{Z}}
\newcommand{\uZ}{\underline{\mathbb{Z}}}
\newcommand{\FF}{\mathscr{F}}
\newcommand{\FP}{\mathscr{P}}
\newcommand{\la}{\langle}
\newcommand{\ra}{\rangle}
\newcommand{\upi}{\underline{\pi}}

\title{The Green functor structure of $RO(D_{2p})$-graded cohomology of a point}
\author{Guoqi Yan}

\begin{document}
\begin{abstract}
    We compute the $RO(D_{2p})$-graded cohomology of a point with constant coefficient $\underline{\mathbb{Z}}$ together with its Green functor structure. Here $D_{2p}$ is the dihedral group with $p$ an odd prime. This result extends the additive computation of Kriz-Lu.
\end{abstract}

\maketitle

\section{Introduction}

The author learned from Kriz-Lu \cite{KrizLu20} that it was a 1982 Northwestern conference problem to compute the $RO(G)$-graded coefficients of a point for a non-abelian group, and that Peter May at a 2017 colloquium at the University of Michigan mentioned that no such computation was know by then.

The purpose of the current work is to illustrate the computational power of the Tate square when dealing with multiplicative structures, as well as finding clean presentation of the structures. Previously, the additive structures of $RO(D_{2p})$-graded cohomology of a point with $\uZ$ and $\underline{\mathbb{A}}$ coefficients were computed in Kriz-Lu \cite{KrizLu20}. The ring structure in the case of $\uZ$-coefficient was also computed by Liu \cite{Liu22} as an application of his splitting method. The author had private conversations with both Lu and Liu and we acknowledge the independence of our computations.

Now we get into the details of our techniques. For any compact Lie group $G$, the following pullback square
\[\xymatrix{
H\underline{\mathbb{Z}}_h\ar[r]\ar[d]_{\simeq} &H\underline{\mathbb{Z}}\ar[r]\ar[d] &\widetilde{H\underline{\mathbb{Z}}}\ar[d]\\
H\underline{\mathbb{Z}}_h\ar[r] &H\underline{\mathbb{Z}}^h\ar[r] &H\underline{\mathbb{Z}}^t
}\]
was introduced by Greenlees-May in \cite{GM95}. Here $\widetilde{H\uZ}=H\uZ\wedge \widetilde{EG},H\uZ^h=F(EG_+,H\uZ),H\uZ_h=EG_+\wedge H\uZ$ and $H\uZ^t=\widetilde{EG}\wedge F(EG_+,H\uZ)$. The second row is computationally accessible via spectral sequences derived from the cellular filtration of $EG_+$. The three spectral sequences are known as the homotopy orbit spectral sequence (HOSS), the homotopy fixed-point spectral sequence (HFPSS) and the Tate spectral sequence respectively.

For the homotopy of $\widetilde{H\uZ}$, we have the following general tool. Recall for any finite group $G$ and any family of subgroups $\FF$, we have the universal space $E\FF$ characterized up to $G$-equivalence by \cite{Alaska}
\[E\FF^H=
\begin{cases}
	*,\quad H\in \FF,\\
	\emptyset,\quad \text{otherwise}.
\end{cases}
\]
The $G$-space $\widetilde{E\FF}$ is defined in the homotopy cofiber sequence
\[
E\FF_+\to S^0\to \widetilde{E\FF}.
\]
Now if we have two families of subgroups $\FF,\FF'$, we get the following homotopy pullback square which is viewed as a diagram in the homotopy category of $G$-spectra via the suspension spectrum functor
\[
\xymatrix{
\widetilde{E(\FF\cap\FF')}\ar[r]\ar[d]&\widetilde{E\FF}\ar[d]\\
\widetilde{E\FF'}\ar[r]&\widetilde{E(\FF\cup\FF')}
}
\]
This can be seen by applying the geometric fixed-point functors $\Phi^H,H\subset G$, using the fact that $\Phi^H$ commute with taking suspension spectra, and they collectively detect weak equivalences of $G$-spectra.

We take the group presentation $G=D_{2p}=\langle\xi,\tau|\xi^p=\tau^2=1,\tau\xi=\xi^{p-1}\tau\rangle$. By \cite[Proposition 4]{KrizLu20}, the $RO(D_{2p})$-grading can be chosen to be
\[
RO(D_{2p})=\Z\{1,\alpha,\gamma\}
\]
without losing essential information. $1$ stands for the $1$-dimentional trivial representation. $\alpha$ stands for the $1$-dimensional sign representation where $\tau$ acts by $-1$ and $\xi$ acts trivially. $\gamma$ denotes the two-dimential irreducible representation where $\tau$ acts by sending $x+iy,x,y\in \mathbb{R}$ to $x-iy$, and $\xi$ acts by multiplication by $e^{2\pi i/p}$.

The main theorem of this paper is
\begin{Thm}\label{main theorem}
The $RO(D_{2p})$-graded coefficients of $H\uZ$ is given by
\begin{align*}
	\pi_{\star}^{G/G}H\underline{\mathbb{Z}}=\mathbb{Z}&[u_{2\alpha},u_{\gamma-\alpha},a_{\alpha},a_{\gamma}]/(2a_{\alpha},pa_{\gamma})\\
	\oplus&\mathbb{Z}[u_{\gamma-\alpha}]\la 2u_{2\alpha}^{-i}\ra _{i\geq 1}\\
	\oplus&\mathbb{Z}[u_{2\alpha}]\la pu_{\gamma-\alpha}^{-i}\ra _{i\geq 1}\\
	\oplus&\mathbb{Z}\la 2pu_{2\alpha}^{-i}u_{\gamma-\alpha}^{-j}\ra _{i,j\geq 1}\\
	\oplus&\mathbb{Z}/p[u_{\gamma-\alpha}]\la u_{2\alpha}^{-j}a_{\gamma}^i\ra _{i,j\geq 1}\\
	\oplus&\mathbb{Z}/2[u_{2\alpha}]\la u_{\gamma-\alpha}^{-j}a_{\alpha}^i\ra _{i,j\geq 1}\\
	\oplus&\mathbb{Z}/p[u_{2\alpha}^{\pm}]\la \Sigma^{-1}u_{\gamma-\alpha}^{-j}a_{\gamma}^{-i}\ra _{i,j\geq 1}\\
	\oplus&\mathbb{Z}/2[u_{\gamma-\alpha}^{\pm}]\la \Sigma^{-1}u_{2\alpha}^{-j}a_{\alpha}^{-i}\ra _{i,j\geq 1}.
\end{align*}
Moreover, the multiplicative and Mackey functor structures are also determined in Section \ref{Ring and Mackey structure}.
\end{Thm}

\subsection*{Notation}
\begin{itemize}
    \item $G=D_{2p}$ with $p$ an odd prime in this paper.
	\item $\langle\quad\rangle$ means additive generators, $[\quad]$ means polynomials. For example, $\mathbb{Z}[a]\la b\ra$ means that we have the torsion free elements $\{a^ib|i\geq 0\}$.;
	\item Otherwise stated, homotopy groups means $G/G$-level;
    \item All the Mackey functors drawn will be in the following shape of lattice of subgroups:
    
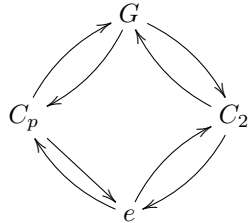
\[\xymatrix{
    	&G\ar@/^/[ld]\ar@/^/[rd]\\
    	C_p\ar@/^/[ur]\ar[dr]&&C_2\ar@/^/[ul]\ar@/^/[dl]\\
    	&e\ar@/^/[ul]\ar@/^/[ur]
    }\]
\end{itemize}

\subsection*{Acknowledgement} The author previously computed the $D_6=\Sigma_3$ case in July 2020. He would like to thank Igor Kriz and Yunze Lu for showing that this implies the general case. He wants to thank Mark Behrens from whom he learned the pullback square of universal spaces of families. He wants to thank for private conversations with Yunze Lu and Yutao Liu. He also benefits from the thesis of Yan Zou to whom he is grateful for.

\section{Euler class and orientation classes}
For a $G$-representation $V$ with $V^G=0$, the inclusion of fixed points $S^0\to S^V$ is not null, which gives us a class $a_V\in \pi_{-V}^{G/G}S^0$. Its Hurewicz image in $H\uZ$ will also be called $a_V$. Theses classes satisfies $a_Va_W=a_{V\oplus W}$.

The other important classes are the orientation classes. Kriz \cite{Sophie} shows that for any finite $G$ and any faithful $G$-representation $V$, we have $\upi_{|V|-V}H\uZ=\uZ$. In our case, we get the orientation classes $u_{\gamma}\in \pi_{2-\gamma}^{G/G}H\uZ$ and $u_{2\alpha}\in\pi^{G/G}_{2-2\alpha}H\uZ$. These classes are also multiplicative, in the sense the $u_{V+W}=u_Vu_W$.

One caveat is that we also have the orientation class $u_{\gamma-\alpha}\in\pi_{1-\gamma+\alpha}^{G/G}H\uZ$ which also generates the Mackey functor $\uZ$. To see this, notice that we can compute $\pi_{1-\gamma+\alpha}^{G/G}H\uZ$ from $\pi_{1-\gamma}^{G/G}H\uZ$ via the cofiber sequence
\[
{G/C_p}_+\to S^0\to S^{\alpha},
\]
and we have
\[
\pi_{1-\gamma}^{G/G}H\uZ=H_1(S^{\gamma};\uZ).
\]
The sphere $S^{\gamma}$ has a cellular structure with a $0$-cell $G/G_+$, two $1$-cells ${G/C_p}_+\wedge S^1$ and one $2$-cell $G/e_+\wedge S^2$. The reduced cellular chain is
\[
0\leftarrow \Z\xleftarrow{\nabla\oplus \nabla}\Z[G/C_p]\oplus \Z[G/C_p]'\xleftarrow{(1,-1)}\Z[G/e].
\]
The map $\nabla$ is the natural augmentation map. The map $(1,-1)$ means sending $x$ to $(\bar{x},-\bar{x})$. Note the first $1$-cell is chosen as the line segment from $0$ to $1$ in the complex plane. We rotate it by $\pi/3$ to get the second $1$-cell. As a result, $\Z[G/C_p]'$ has a different $\Z[G]$-module structure. More explicitly, we have $\xi\cdot \bar{1}=\bar{\xi}$ and $\tau\cdot\bar{1}=\bar{\xi}^{p-1}$.

Taking $(-)^G$ and then homology, we get $H_0(S^{\gamma})=\Z/3,H_1(S^{\gamma})=\Z/2$ and $H_2(S^{\gamma})=0$. We look at the long exact sequence in homology derived from the cofiber sequence
\[
S^{\gamma-\alpha}\to S^{\gamma}\to {G/C_p}_+\wedge S^{\gamma}.
\]
We deduce that we have the short exact sequence
\[
0\to \Z\to \pi^{G/G}_{1-\gamma+\alpha}H\uZ\to \Z/2\to 0.
\]
Suppose it splits, and $\pi^{G/G}_{1-\gamma+\alpha}H\uZ=\Z\oplus\Z/2$. Let $x$ be the generator of $\Z/2$. Then $Res^G_{C_2}(x)\in \pi^{G/C_2}_{1-\gamma+\alpha}H\uZ\cong \pi^{G/C_2}_{0}H\uZ\cong \Z$, which implies $Res^G_{C_2}(x)=0$ since $2x=0$. This violates the cohomological condition, since $Tr^G_{C_2}Res^G_{C_2}(x)=px=x$. We conclude $\pi^{G/G}_{1-\gamma+\alpha}H\uZ=\Z$. 

We can easily show that the Mackey functor structure is $\uZ$. We already know $Res^{C_p}_e=Res^{C_2}_e=1$. By the cohomological condition, we know $Res^G_{C_p}=1$ or $2$ and $Res^G_{C_p}=1$ or $p$. By the transitivity of restrictions, the only option is $Res^G_{C_p}=Res^G_{C_2}=1$.

\section{Computations in the second row}
\subsection{The homotopy of \texorpdfstring{$H\uZ^h$}{} via the HFPSS}
The HFPSS takes the form
\[
E_2^{s,V}=H^s(G;\pi^{G/e}_V H\uZ)\Rightarrow \pi^{G/G}_{V-s}H\uZ^h.
\]
It collapses since $\pi^{G/e}_V H\uZ$ concentrates at gradings with virtual dimension $0$. To understand the $E_2$-page, we look at the Lyndon-Hochschild-Serre spectral sequence. We have
\[
E_2^{i,j}=H^i(G/{C_p};H^j(C_p;M))\Rightarrow H^{i+j}(G;M)
\]
for $M$ a $G$-module. $G/{C_p}\cong C_2\la \overline{\tau}\ra $.

If $M=\mathbb{Z}$, then $\overline{\tau}$ act as $-1$ on $H^2(C_p;\mathbb{Z})$, and act as ring homomorphisms since $\mathbb{Z}$ is a ring object in the category of $\mathbb{Z}[G]$-modules.

If $M=\widetilde{\mathbb{Z}}$ (the sign representation of $G$), then $\overline{\tau}$ act as $-1$ on $H^0(C_p;\widetilde{\mathbb{Z}}),H^4(C_p;\widetilde{\mathbb{Z}})$ and trivially on $H^2(C_p;\widetilde{\mathbb{Z}})$, and it is $4$-periodic. Note $G/{C_p}$ does not act as ring homomorphisms since $\widetilde{\mathbb{Z}}$ is not a ring object in the category of $\mathbb{Z}[G]$-modules. But anyway, $H^*(G;\widetilde{\mathbb{Z}})$ is a module over $H^*(G;\mathbb{Z})$. We have
\begin{align*}
&H^*(G;\mathbb{Z})=\mathbb{Z}[x,z]/{(2x,pz)},|x|=2,|z|=4;\\
&H^*(G;\widetilde{\mathbb{Z}})=\mathbb{Z}[x]/{(2x)}\la c\ra \oplus \mathbb{Z}/p[z]/{(pz)}\la d\ra ,|c|=1,|d|=2.
\end{align*}
As Mackey functors, they are:
\[
\underline{H}^*(G;\mathbb{Z})
\]
\[\xymatrix{
	&\mathbb{Z}[x,z]/(2x,pz)\ar@/^/[ld]^{x\mapsto 0,z\mapsto y^2}\ar@/^/[rd]^{x\mapsto w,z\mapsto 0}\\
	\mathbb{Z}[y]/py,|y|=2,\bar{\tau}\cdot y=-y\ar@/^/[ur]^{y\mapsto 0,y^2\mapsto 2z}\ar[dr]&&\mathbb{Z}[w]/2w,|w|=2\ar@/^/[ul]^{w\mapsto x}\ar@/^/[dl]\\
	&\mathbb{Z}\ar@/^/[ul]\ar@/^/[ur]
}\]
Note $Res^G_{C_p}$ identifies $H^*(G;\mathbb{Z})_{(p)}$ (the p-local part of $H^*(G;\mathbb{Z})$) as the subring $\mathbb{Z}[y^2]/py^2$ of $H^*(C_p;\mathbb{Z})$.

\[
\underline{H}^*(G;\widetilde{\mathbb{Z}})
\]
\[\xymatrix{
	&H^*(G;\widetilde{\mathbb{Z}})=\mathbb{Z}[x]/{(2x)}\la c\ra \oplus \mathbb{Z}/p[z]/{(pz)}\la d\ra \ar@/^/[ld]^{x\mapsto 0,z\mapsto y^2,d\mapsto y}\ar@/^/[rd]^{x\mapsto w,z\mapsto 0,c\mapsto e}\\
	\mathbb{Z}[y]/py\ar@/^/[ur]^{y\mapsto 0,y^2\mapsto 2z}\ar[dr]&&\mathbb{Z}[w]/2w\la e\ra ,|e|=2\ar@/^/[ul]^{w\mapsto x,e\mapsto c}\ar@/^/[dl]\\
	&\mathbb{Z}\ar@/^/[ul]\ar@/^/[ur]
}\]
Note $Res^G_{C_p}$ identifies $H^*(G;\widetilde{\mathbb{Z}})_{(p)}$ as the subring $\mathbb{Z}[y^2]/py^2\la y\ra $ of $H^*(C_p;\widetilde{\mathbb{Z}})$. Note the structure maps above satisfy both the cohomological condition and the double coset formula.

From the HFPSS, we get
\[
E_2^{s,t,i,j}=H^{-s}(G;\pi_{t+i\alpha+j\gamma}^{G/e}(H\underline{\mathbb{Z}}))\Rightarrow \pi_{s+t+i\alpha+j\gamma}^{G/G}(H\underline{\mathbb{Z}}^h);
\]
We conclude
\begin{equation}\label{H^h}
H\underline{\mathbb{Z}}^h_{\star}=\mathbb{Z}[a_{\alpha},a_{\gamma},u_{2\alpha}^{\pm},u_{\gamma-\alpha}^{\pm}]/{(2a_{\alpha},pa_{\gamma})}.    
\end{equation}

\subsection{The rest of the computations in the second row}

As mentioned in the introduction, we have the following pullback squares in $Ho(Sp^G)$:
\[\xymatrix{
\widetilde{EG}\ar[d]\ar[r] &\widetilde{E\{e,C_p\}}\ar[d]\\
\widetilde{E\{e,C_2\}}\ar[r] &\widetilde{E\FP}
}\]
Since $\widetilde{E\FP}\simeq S^{\infty\alpha}\wedge S^{\infty\gamma}$, we conclude $H\underline{\mathbb{Z}}\wedge\widetilde{E\FP}\simeq *$, because $H\underline{\mathbb{Z}}\wedge\widetilde{E\FP}\simeq H\underline{\mathbb{Z}}[a_{\alpha}^{-1},a_{\gamma}^{-1}]$, and $a_{\alpha},a_{\gamma}$ are 2-torsion and p-torsion respectively. Thus we have
\begin{equation}\label{Splitting of tildeHZ}
\widetilde{H\uZ}\simeq (H\underline{\mathbb{Z}}\wedge\widetilde{E\{e,C_p\}})\vee (H\underline{\mathbb{Z}}\wedge\widetilde{E\{e,C_2\}})\simeq H\underline{\mathbb{Z}}[a_{\alpha}^{-1}]\vee H\underline{\mathbb{Z}}[a_{\gamma}^{-1}].
\end{equation}

The same splitting happens when we replace $H\uZ$ by $H\uZ^h$. From \eqref{H^h}, we conclude 
\begin{align*}
	H\underline{\mathbb{Z}}^t_{\star}=&\mathbb{Z}/2[a_\alpha^{\pm},u_{2\alpha}^{\pm},u_{\gamma-\alpha}^{\pm}]\\
	&\oplus\mathbb{Z}/p[a_{\gamma}^{\pm},u_{2\alpha}^{\pm},u_{\gamma-\alpha}^{\pm}].
\end{align*}
In the following sequence, let $K,C$ denote kernel and cokernel respectively,
\[
0\to K\to H\underline{\mathbb{Z}}^h_{\star}\to H\underline{\mathbb{Z}}^t_{\star}\to C\to 0.
\]
We have
\begin{align*}
	K=&2p\mathbb{Z}[u_{2\alpha}^{\pm},u_{\gamma-\alpha}^{\pm}]\\
	C=&\mathbb{Z}/2[u_{2\alpha}^{\pm},u_{\gamma-\alpha}^{\pm}]\la a_{\alpha}^{-i}\ra _{i\geq 1}\\
	&\oplus\mathbb{Z}/p[u_{2\alpha}^{\pm},u_{\gamma-\alpha}^{\pm}]\la a_{\gamma}^{-i}\ra _{i\geq 1}
\end{align*}
There is no extension problem in the following sequence since $K$ is free.
\[
0\to \Sigma^{-1}C\to H\underline{\mathbb{Z}}_{h\star}\to K\to 0
\]
Finally, we conclude
\begin{align*}
	H\underline{\mathbb{Z}}_{h\star}=&2p\mathbb{Z}[u_{2\alpha}^{\pm},u_{\gamma-\alpha}^{\pm}]\\
	&\oplus\mathbb{Z}/2[u_{2\alpha}^{\pm},u_{\gamma-\alpha}^{\pm}]\la \Sigma^{-1}a_{\alpha}^{-i}\ra _{i\geq 1}\\
	&\oplus\mathbb{Z}/p[u_{2\alpha}^{\pm},u_{\gamma-\alpha}^{\pm}]\la \Sigma^{-1}a_{\gamma}^{-i}\ra _{i\geq 1}
\end{align*}

\section{Computations in the first row}
Next we compute $\widetilde{H\underline{\mathbb{Z}}}_{\star}$ using the splitting \eqref{Splitting of tildeHZ}.

\subsection{Computing \texorpdfstring{$H\underline{\mathbb{Z}}[a_{\alpha}^{-1}]_{\star}$}{}}
From the fact that
\[
H\underline{\mathbb{Z}}[a_{\alpha}^{-1}]_{a+b\alpha}^{G/C_2}=[{G/C_2}_+\wedge S^{a+b\alpha},H\underline{\mathbb{Z}}\wedge S^{\infty\alpha}]^G=[S^{a+b\alpha},H\underline{\mathbb{Z}}\wedge S^{\infty\alpha}]^{C_2}
\]
we get
\[
H\underline{\mathbb{Z}}[a_{\alpha}^{-1}]^{G/C_2}_{*+*\alpha}=\mathbb{Z}/2[a_{\alpha}^{\pm},u_{2\alpha}].
\]
We notice that the Mackey functor $H\underline{\mathbb{Z}}[a_{\alpha}^{-1}]_{\star}$ is 2-torsion on each level. Since $H\underline{\mathbb{Z}}[a_{\alpha}^{-1}]_{\star}$ is obviously an $H\underline{\mathbb{Z}}$-module, meaning that its homotopy groups is cohomological, we know $Tr^G_{C_2}Res^G_{C_2}=p$ is an isomorphism on the $G/G$-level when it is not 0. By the double coset formula, we get $Res^G_{C_2}Tr^G_{C_2}=p$ is an isomorphism on the $G/{C_2}$-level when it is not 0. We conclude that
\[
H\underline{\mathbb{Z}}[a_{\alpha}^{-1}]^{G/G}_{*+*\alpha}=\mathbb{Z}/2[a_{\alpha}^{\pm},u_{2\alpha}]
\]
and the Mackey functor structure is:
\[\xymatrix{
&\mathbb{Z}/2[a_{\alpha}^{\pm},u_{2\alpha}]\ar@/^/[ld]\ar@/^/[rd]^1\\
0\ar@/^/[ur]\ar[dr]&&\mathbb{Z}/2[a_{\alpha}^{\pm},u_{2\alpha}]\ar@/^/[ul]^1\ar@/^/[dl]\\
&0\ar@/^/[ul]\ar@/^/[ur]
}\]
Note that in the same way we can deduce $Res^G_{C_2}=Tr^G_{C_2}=1$ for the entire gradings $\star$, not just for ${*+*\alpha}$.\\

\textbf{Remark:} There is another way to know $H\underline{\mathbb{Z}}[a_{\alpha}^{-1}]_{*}$ (thus $H\underline{\mathbb{Z}}[a_{\alpha}^{-1}]_{*+*\alpha}$ by $a_{\alpha}$-periodicity). We can first compute 
\begin{align*}
	\widetilde{H\underline{\mathbb{Z}}}_*=&\mathbb{Z}/p\la (\frac{u_{2\gamma}}{a_{2\gamma}})^k\ra _{k\geq 1}\\
	&\oplus \mathbb{Z}/2\la (\frac{u_{2\alpha}}{a_{2\alpha}})^i\ra _{i\geq 1}\\
	&\oplus \mathbb{Z}/2p\la 1\ra .
\end{align*}

The splitting $\widetilde{H\underline{\mathbb{Z}}}_*\simeq H\underline{\mathbb{Z}}[a_{\alpha}^{-1}]_{*}\vee H\underline{\mathbb{Z}}[a_{\gamma}^{-1}]_{*}$ is a splitting into 2-torsion and p-torsion. So we must have $\mathbb{Z}/2p=\mathbb{Z}/p\oplus \mathbb{Z}/2$, and
\begin{align*}
H\underline{\mathbb{Z}}[a_{\alpha}^{-1}]_{*}=&\mathbb{Z}/2[\frac{u_{2\alpha}}{a_{2\alpha}}]\\
H\underline{\mathbb{Z}}[a_{\gamma}^{-1}]_{*}=&\mathbb{Z}/p[\frac{u_{2\gamma}}{a_{2\gamma}}].
\end{align*}

Now we want to understand the $RO(G)$-gradings which involves $\gamma$. Note that under the restriction map $i^*:RO(G)\to RO(C_2)$, we have $i^*(\gamma)=1+\alpha$. That means the restriction $Res^G_{C_2}:\pi^{G/G}_{1-\gamma+\alpha}H\underline{\mathbb{Z}}[a_{\alpha}^{-1}]\to \pi^{G/C_2}_{1-\gamma+\alpha}H\underline{\mathbb{Z}}[a_{\alpha}^{-1}]\cong \pi^{G/G}_{0}H\underline{\mathbb{Z}}[a_{\alpha}^{-1}]$ send $u_{\gamma-\alpha}$ to $1$. On the $G/C_2$-level, $u_{\gamma-\alpha}$ is invertible so that it has an inverse $u_{\gamma-\alpha}^{-1}\in\pi^{G/C_2}_{-1+\gamma-\alpha}$, which lifts to the $G/G$-level along the restriction map. If we call this lift $\tilde{u}_{\gamma-\alpha}$, we have $Res^G_{C_2}(\tilde{u}_{\gamma-\alpha}\cdot u_{\gamma-\alpha})=Res^G_{C_2}(\tilde{u}_{\gamma-\alpha})\cdot Res^G_{C_2}(u_{\gamma-\alpha})=1$. Since restriction is an isomorphism, we conclude $\tilde{u}_{\gamma-\alpha}\cdot u_{\gamma-\alpha}=1$ thus $\tilde{u}_{\gamma-\alpha}$ deserves the name $u_{\gamma-\alpha}^{-1}$. Hence $H\underline{\mathbb{Z}}[a_{\alpha}^{-1}]_{\star}$ is $u_{\gamma-\alpha}$-local and its $RO(G)$-graded homotopy Mackey functor $\underline{H\underline{\mathbb{Z}}[a_{\alpha}^{-1}]_{\star}}$ is 

\[\xymatrix{
	&\mathbb{Z}/2[a_{\alpha}^{\pm},u_{2\alpha},u_{\gamma-\alpha}^{\pm}]\ar@/^/[ld]\ar@/^/[rd]^1\\
	0\ar@/^/[ur]\ar[dr]&&\mathbb{Z}/2[a_{\alpha}^{\pm},u_{2\alpha}]\ar@/^/[ul]^1\ar@/^/[dl]\\
	&0\ar@/^/[ul]\ar@/^/[ur]
}\]

\subsection{Computing \texorpdfstring{$H\underline{\mathbb{Z}}[a_{\gamma}^{-1}]_{\star}$}{}}
We start from
\[
H\underline{\mathbb{Z}}[a_{\gamma}^{-1}]_{*}=\mathbb{Z}/p[\frac{u_{2\gamma}}{a_{2\gamma}}]
\]
By $a_{\gamma}$-periodicity, we know $H\underline{\mathbb{Z}}[a_{\gamma}^{-1}]_{4*+*\gamma}=\mathbb{Z}/p[a_{\gamma}^{\pm},u_{2\gamma}],*\in\mathbb{Z}$. We claim that this spectrum is $u_{2\alpha}$-local. By the cohomological condition and double coset formula, we have
\[
\pi_{\pm(2-2\alpha)}^{G/G}H\underline{\mathbb{Z}}[a_{\gamma}^{-1}]\cong\pi_{\pm(2-2\alpha)}^{G/C_p}H\underline{\mathbb{Z}}[a_{\gamma}^{-1}]\cong\pi_{0}^{G/C_p}H\underline{\mathbb{Z}}[a_{\gamma}^{-1}]\cong\Z/p.
\]
We also have $Res^G_{C_p}(u_{2\alpha})$ is a unit. Let the lift of its inverse to the $G/G$-level be $\tilde{u_{2\alpha}}$. As before, we can argue that $\tilde{u_{2\alpha}}\cdot u_{2\alpha}=1$, and $\tilde{u_{2\alpha}}$ deserves the name $u_{2\alpha}^{-1}$. As a result, it suffices to understand the gradings 
\[
\underline{H\underline{\mathbb{Z}}[a_{\gamma}^{-1}]_{*+*\gamma}}\text{ 
 and }\underline{H\underline{\mathbb{Z}}[a_{\gamma}^{-1}]_{*+*\gamma+\alpha}}.
 \]

The Mackey functor structure of
$\underline{H\underline{\mathbb{Z}}[a_{\gamma}^{-1}]_{*+*\gamma}}$ is
\[\xymatrix{
	&\mathbb{Z}/p[a_{\gamma}^{\pm},u_{2\gamma}]\ar@/^/[ld]^-{1}\ar@/^/[rd]\\
	\mathbb{Z}/p[a_{\gamma}^{\pm},u_{\gamma}]\ar@/^/[ur]^-{2(0)}\ar[dr]&&0\ar@/^/[ul]\ar@/^/[dl]\\
	&0\ar@/^/[ul]\ar@/^/[ur]
}\]
Here the map $2(0)$ sends monomials with odd powers of $u_{\lambda}$ to $0$ and monomials with even powers of $u_{\lambda}$ to $2$ times the corresponding monomial on the $G/G$-level.

Now since we have $Res^G_{C_p}(u_{\gamma-\alpha})=u_{\gamma}\in\pi^{G/C_p}_{2-\gamma}$, we know $Res^G_{C_p}:\mathbb{Z}/p[u_{2\gamma}]\la u_{\gamma-\alpha}\ra \mapsto \mathbb{Z}/p[u_{\gamma}^2]\la u_{\gamma}\ra $. So the Mackey functor $\underline{H\underline{\mathbb{Z}}[a_{\gamma}^{-1}]_{*+*\gamma+\alpha}}$ is 
\[\xymatrix{
	&\mathbb{Z}/p[a_{\gamma}^{\pm},u_{2\gamma}]\la 1,u_{\gamma-\alpha}\ra \ar@/^/[ld]^-{1}\ar@/^/[rd]\\
	\mathbb{Z}/p[a_{\gamma}^{\pm},u_{\gamma}]\ar@/^/[ur]^-{2(0)}\ar[dr]&&0\ar@/^/[ul]\ar@/^/[dl]\\
	&0\ar@/^/[ul]\ar@/^/[ur]
}\]

Using the relation $u_{2\alpha}\cdot(u_{\gamma-\alpha})^2=u_{2\gamma}$, we have $\underline{H\underline{\mathbb{Z}}[a_{\gamma}^{-1}]}_{\star}=\mathbb{Z}/p[a_{\gamma}^{\pm},u_{\gamma-\alpha},u_{2\alpha}^{\pm}]$, and the entire Mackey functor $\underline{H\underline{\mathbb{Z}}[a_{\gamma}^{-1}]_{\star}}$ is
\[\xymatrix{
	&\mathbb{Z}/p[a_{\gamma}^{\pm},u_{\gamma-\alpha},u_{2\alpha}^{\pm}]\ar@/^/[ld]^-{1}\ar@/^/[rd]\\
	\mathbb{Z}/p[a_{\gamma}^{\pm},u_{\gamma}]\ar@/^/[ur]^-{2(0)}\ar[dr]&&0\ar@/^/[ul]\ar@/^/[dl]\\
	&0\ar@/^/[ul]\ar@/^/[ur]
}\]

\section{The homotopy of \texorpdfstring{$H\underline{\mathbb{Z}}$}{}}
Now we have
\begin{align*}
	\widetilde{H\underline{\mathbb{Z}}}_{\star}=&\mathbb{Z}/2[u_{2\alpha},a_{\alpha}^{\pm},u_{\gamma-\alpha}^{\pm}]\\
	&\oplus\mathbb{Z}/p[u_{2\alpha}^{\pm},a_{\gamma}^{\pm},u_{\gamma-\alpha}]
\end{align*}
Taking kernels and cokernels of the following
\[
0\to K'\to \widetilde{H\underline{\mathbb{Z}}}_{\star}\to H\underline{\mathbb{Z}}_{h\star}\to C'\to 0
\]
we get
\begin{align*}
	K'=&\mathbb{Z}/p[{a_{\gamma},u_{2\alpha}^{\pm},u_{\gamma-\alpha}}]\\
	&\oplus\mathbb{Z}/2[a_{\alpha},u_{2\alpha},u_{\gamma-\alpha}^{\pm}]
\end{align*}
and
\begin{align*}
	C'=&2p\mathbb{Z}[u_{2\alpha}^{\pm},u_{\gamma-\alpha}^{\pm}]\\
	&\oplus\mathbb{Z}/p[u_{2\alpha}^{\pm}]\la \Sigma^{-1}u_{\gamma-\alpha}^{-j}a_{\gamma}^{-i}\ra _{i,j\geq 1}\\
	&\oplus\mathbb{Z}/2[u_{\gamma-\alpha}^{\pm}]\la \Sigma^{-1}u_{2\alpha}^{-j}a_{\alpha}^{-i}\ra _{i,j\geq 1}.
\end{align*}

Now we have the extension problem
\[
0\to C'\to H\underline{\mathbb{Z}}_{\star}\to K'\to 0
\]
For degree reasons, we do not need to consider the elements with a $\Sigma^{-1}$. The elements in $K'$ with positive power of $a$'s contribute directly to $H\underline{\mathbb{Z}}_{\star}$. 

Then we consider the non-trivial extensions: We have four different cases depending on the number of $u$'s being inverted. Let $A=\pi_V^{G/G}(H\underline{\mathbb{Z}})$ be the unknown group in the corresponding dimension.

(1) We have
\[
0\to 2p\mathbb{Z}[u_{2\alpha},u_{\gamma-\alpha}]\to A\to \mathbb{Z}/2p[u_{2\alpha},u_{\gamma-\alpha}]\to 0.
\]
Since this is an extension of $H\underline{\mathbb{Z}}_{\star}$-module, we only need to consider
\[
0\to 2p\mathbb{Z}\la 1\ra \to A\to \mathbb{Z}/2p\la 1\ra \to 0. 
\]
$A=\mathbb{Z}$ since this is just $H\underline{\mathbb{Z}}_0$.

(2) We have
\[
0\to 2p\mathbb{Z}[u_{\gamma-\alpha}]\la u^{-i}_{2\alpha}\ra _{i\geq 1}\to A\to \mathbb{Z}/p[u_{\gamma-\alpha}]\la u^{-i}_{2\alpha}\ra _{i\geq 1}\to 0
\]
Similarly, we consider
\[
0\to 2p\mathbb{Z}\la u^{-i}_{2\alpha}\ra _{i\geq 1}\to A\to \mathbb{Z}/p\la u^{-i}_{2\alpha}\ra _{i\geq 1}\to 0
\]
The group is in degree $i(2\alpha-2)$, and both $H\underline{\mathbb{Z}}_{i(2\alpha-2)}$ and $H\underline{\mathbb{Z}}^h_{i(2\alpha-2)}$ are $\underline{\mathbb{Z}}$. The underlying map from the former to the latter is $1$. Thus in this degree, the map is $1:\underline{\mathbb{Z}}\to \underline{\mathbb{Z}}$. Now we know the extension is
\[
0\to 2p\mathbb{Z}\la u^{-i}_{2\alpha}\ra _{i\geq 1}\xrightarrow{p} \mathbb{Z}\la x\ra \to \mathbb{Z}/p\la u^{-i}_{2\alpha}\ra _{i\geq 1}\to 0
\]
Now $2pu^{-i}_{2\alpha}$ maps to $2p\cdot u^{-i}_{2\alpha}=p\cdot x\in H\underline{\mathbb{Z}}^h_{i(2\alpha-2)}$, we get $x=2u^{-i}_{2\alpha}$.

(3) We have
\[
0\to 2p\mathbb{Z}[u_{2\alpha}]\la u^{-i}_{\gamma-\alpha}\ra _{i\geq 1}\to A\to \mathbb{Z}/2[u_{2\alpha}]\la u^{-i}_{\gamma-\alpha}\ra _{i\geq 1}\to 0
\]
Similarly, we consider
\[
0\to 2p\mathbb{Z}\la u^{-i}_{\gamma-\alpha}\ra _{i\geq 1}\to A\to \mathbb{Z}/2\la u^{-i}_{\gamma-\alpha}\ra _{i\geq 1}\to 0
\]
The group is in degree $i(-1+\gamma-\alpha)$, and both $H\underline{\mathbb{Z}}_{i(-1+\gamma-\alpha)}$ and $H\underline{\mathbb{Z}}^h_{i(-1+\gamma-\alpha)}$ are $\underline{\mathbb{Z}}$. In this degree, the map is $1:\underline{\mathbb{Z}}\to \underline{\mathbb{Z}}$. Now we know the extension is
\[
0\to 2p\mathbb{Z}\la u^{-i}_{\gamma-\alpha}\ra _{i\geq 1}\xrightarrow{2} \mathbb{Z}\la y\ra \to \mathbb{Z}/2\la u^{-i}_{\gamma-\alpha}\ra _{i\geq 1}\to 0
\]
Now $2pu^{-i}_{\gamma-\alpha}$ maps to $2p\cdot u^{-i}_{\gamma-\alpha}=2\cdot y\in H\underline{\mathbb{Z}}^h_{i(-1+\gamma-\alpha)}$, we get $y=pu^{-i}_{\gamma-\alpha}$.

(4) We have
\[
0\to 2p\mathbb{Z}\la u_{2\alpha}^{-j}u^{-i}_{\gamma-\alpha}\ra _{i,j\geq 1}\to A\to 0
\]
from which we get $A=\mathbb{Z}\la 2pu_{2\alpha}^{-j}u^{-i}_{\gamma-\alpha}\ra _{i,j\geq 1}$.

As a conclusion, we proved Theorem \ref{main theorem}.

\section{The multipicative structure and Mackey functor structure\label{Ring and Mackey structure}}
\subsection{Multiplicative structure}
Since every spectrum in the Tate square is an $H\uZ$-module and every map is a $H\uZ$-module map, the multiplications by Euler and orientation classes are obvious in Theorem \ref{main theorem}. 

For torsion-free classes, the number in front of them should be explained as follows. Take classes in $\Z[u_{\gamma-\alpha}]\la 2u_{2\alpha}^{-i}\ra _{i\geq 1}$ as an example. The $2$ means that the images of these classes under either the map $H\uZ\to H\uZ^h$ or $Res^G_e$ is $2$ times of an actual torsion-free class. Similarly for other torsion-free classes. So the product of any two of such classes are determined by their images under Borel completion or $Res^G_e$. For example
\[
(2u_{\gamma-\alpha}u_{2\alpha}^{-1})^2=2\cdot 2u_{\gamma-\alpha}^2u_{2\alpha}^{-2}\text{ and }2u_{\gamma-\alpha}u_{2\alpha}^{-1}\cdot pu_{\gamma-\alpha}^{-1}=p\cdot 2u_{2\alpha}^{-1}. 
\]

For torsion classes, firstly we have that the product of a $2$-torsion class and a $p$-torsion class is always $0$. Then we have that the product of any two classes with a $\Sigma^{-1}$ in front of them is $0$. Take the product of any two classes from $\mathbb{Z}/p[u_{2\alpha}^{\pm}]\la \Sigma^{-1}u_{\gamma-\alpha}^{-j}a_{\gamma}^{-i}\ra _{i,j\geq 1}$ as an example. The product should be $p$-torsion and killed by powers of $a_{\lambda}$, so the only possibility is that it still belongs to $\mathbb{Z}/p[u_{2\alpha}^{\pm}]\la \Sigma^{-1}u_{\gamma-\alpha}^{-j}a_{\gamma}^{-i}\ra _{i,j\geq 1}$. For degree reasons, that cannot happen.

For the rest of the products, they follow the rule of multiplication by names of elements and degree counting. For example
\[
u_{2\alpha}^{-1}a_{\gamma}\cdot \Sigma^{-1}u_{\gamma-\alpha}^{-1}a_{\gamma}^{-2}=u_{2\alpha}^{-1}\Sigma^{-1}u_{\gamma-\alpha}^{-1}a_{\gamma}^{-1}\text{ and }u_{2\alpha}^{-1}a_{\gamma}\cdot \Sigma^{-1}u_{\gamma-\alpha}^{-1}a_{\gamma}^{-1}=0.
\]

\subsection{Mackey functor structure}
The Mackey functor structure of $\upi_{\star}H\uZ$ for the group $C_p$ ($p$ any prime) is well known, see for example \cite{Zeng17}. For the $G$-Mackey functor, we keep the following picture in mind
\[\xymatrix{
	&\pi_{\star}^{G/G}H\uZ \ar@/^/[ld]^-{}\ar@/^/[rd]\\
	\pi_{\star}^{G/C_p}H\uZ[u_{\alpha}^{-1}]\ar@/^/[ur]^-{}\ar[dr]&&\pi_{\star}^{G/C_2}H\uZ[u_{\gamma-\alpha}^{-1}]\ar@/^/[ul]\ar@/^/[dl]\\
	&\pi_{\star}^{G/e}H\uZ[u_{\alpha}^{-1},u_{\gamma-\alpha}^{-1}]\ar@/^/[ul]\ar@/^/[ur].
}\]
Since Euler classes and orientation classes restrict to Euler classes and orientation classes respectively, the Mackey functors generated by $\mathbb{Z}[u_{2\alpha},u_{\gamma-\alpha},a_{\alpha},a_{\gamma}]/(2a_{\alpha},pa_{\gamma})$ are obvious.

For torsion-free classes, take $\mathbb{Z}\la 2pu_{2\alpha}^{-i}u_{\gamma-\alpha}^{-j}\ra _{i,j\geq 1}$ as an example. The corresponding classes at the $G/C_p$-level are $\mathbb{Z}\la pu_{\gamma-\alpha}^{-j}\ra _{j\geq 1}\la u_{2\alpha}^{-i}\ra_{i\geq 1}$. By looking at their restrictions to the $G/e$-level, we know $Res^G_{C_p}=2$. Similarly we deduce $Res^G_{C_2}=p$.

For torsion classes, take $\mathbb{Z}/p[u_{2\alpha}^{\pm}]\la \Sigma^{-1}u_{\gamma-\alpha}^{-j}a_{\gamma}^{-i}\ra _{i,j\geq 1}$ as an example. The corresponding classes at the $G/C_p$-level are $\mathbb{Z}/p[u_{2\alpha}^{\pm}]\la \Sigma^{-1}u_{\gamma}^{-j}a_{\gamma}^{-i}\ra _{i,j\geq 1}$ and we easily deduce $Res^G_{C_p}=1$. By cohomological condition, we know $Tr^G_{C_p}=2$. The classes at other levels are $0$ and there is no non-trivial restriction or transfer to consider.

\bibliographystyle{alpha}
\bibliography{bib}
\end{document}